\documentclass{amsart}
\usepackage{amssymb}
\usepackage{latexsym}
\usepackage{color}
\def\ekaterina{}
\def\ivacomments{}
\begin{document}
\newtheorem{theorem}{Theorem}[section]
\newtheorem{lemma}[theorem]{Lemma}
\newtheorem{proposition}[theorem]{Proposition}
\newtheorem{remark}[theorem]{Remark}
\newtheorem{corollary}[theorem]{Corollary}
\newtheorem{extlemma}[theorem]{Extension Lemma}
\font\pbglie=eufm10
\def\Tr{\operatorname{Tr}}
\def\id{\operatorname{Id}}
\def\ffrac#1#2{{\textstyle\frac{#1}{#2}}}
\def\la{\langle}
\def\ra{\rangle}

\def\ll{\lambda}

\makeatletter
 \renewcommand{\theequation}{%
 \thesection.\alph{equation}}
 \@addtoreset{equation}{section}
 \makeatother
\title[Curvature structure of Self-dual 4-manifolds]
{Curvature structure of Self-dual 4-manifolds}
\author{Novica  Bla{\v z}i{\'c}, Peter Gilkey, Stana Nik\v cevi\'c, and Iva
Stavrov}
\begin{address}{PG: Mathematics Department, University of Oregon,
Eugene Or 97403 USA.\newline Email: {\it gilkey@darkwing.uoregon.edu}}
\end{address}
\begin{address}
{SN: Mathematical Institute, Sanu,
Knez Mihailova 35, p.p. 367,
11001 Belgrade,
Serbia. Email: {\it stanan@mi.sanu.ac.yu}}\end{address}
\begin{address}{IS: 
Department of Mathematical Sciences, Lewis and Clark College, 0615 SW Palatine Hill Road, MSC 110, Portland, Oregon, 97219 USA.
Email: {\it istavrov@lclark.edu}}
\end{address}


\begin{abstract}
 We show the existence of a modified $\mathrm{Cliff}(1,1)$-structure compatible with an Osserman $0$-model of signature $(2,2)$.
We then apply this algebraic result to certain classes of pseudo-Riemannian manifolds of signature $(2,2)$. We obtain a new
characterization of the Weyl curvature tensor of an (anti-)self-dual manifold and we prove some new results regarding (Jordan)
Osserman manifolds.
\end{abstract}

\maketitle
 
\section*{Dedication}
{\ekaterina This paper is one of several projects that were begun by Novica Bla{\v z}i{\' c} but 
{\ivacomments not} 
completed owing to his
untimely death in 2005. The work has been finished to preserve his mathematical legacy and is dedicated to his
memory.}

\section{Introduction}\label{intro}
Let $R$ be the curvature operator of the Levi-Civita connection of a pseudo-Riemannian manifold $M$. The {\it Jacobi operator}
$J_R(x):T_PM\rightarrow T_PM$ corresponding to unit spacelike or unit timelike tangent vectors
$x$ is characterized by
$$J_R(x)y:=R(y,x)x,\ \ x\in T_PM.$$
It plays a central role in curvature theory \cite{ X3,X2}.

It was conjectured by Osserman \cite{Oss} in the Riemannian setting that the spectrum of the Jacobi operator $J_R(x)$ is independent of
the choice of a unit tangent vector
$x$ and its base point $P$ if and only if the underlying Riemannian manifold is locally {\ekaterina rank $1$-symmetric} or flat. This
conjecture follows from work of Chi \cite{C} and Nikolayevsky \cite{Nik0, Nik1, Nik2} in dimensions other than $16$; {\ekaterina this
question is still open in dimension
$16$}. Nikolayevsky executed the approach outlined in \cite{GSV}, the major part of which is in showing that an Osserman
$0$-model (see the next section for the definition) allows a compatible Clifford algebra structure. This crucial step is the
algebraic counterpart to studying the so-called {\it point-wise Osserman manifolds}: manifolds where the spectrum of the Jacobi
operator
$J_R(x)$ does not depend on the choice of the unit tangent vector $x$ at any of the base points, but is allowed to vary from
point to point. Note that there exists a nice connection between point-wise Osserman and self-dual Einstein Riemannian manifolds
due to the work of Sekigawa and Vanhecke \cite{SV}. There are many other properties of the curvature operator which can be studied
similarly -- see, for example, \cite{X1,X5}.

Pseudo Riemannian and conformal geometry is central to many investigations and the phenomena are often very different from the
Riemannian setting (see, for example,
\cite{X7,X8}), and the study of Osserman manifolds is no exception. In the pseudo-Riemannian geometry the study of Osserman-type
manifolds becomes rather complicated because their Jacobi operators need not be diagonalizable. We say a pseudo-Riemannian manifold is
(point-wise) timelike/spacelike Jordan Osserman if the Jordan normal form of $J_R(x)$ is independent of the choice of unit
timelike/spacelike vector
$x$. Numerous examples have been constructed
\cite{BGV} which show the existence of non-homogeneous point-wise Jordan Osserman manifolds. 

The Weyl curvature tensor $W$, which depends only on the conformal class of a 
{\ivacomments pseudo-Riemannian} 
manifold, obeys the same algebraic symmetries as
the Riemann curvature tensor itself (see the next section for details). Quite naturally one is led to investigating the
spectral geometry of the {\it conformal Jacobi operator}
$$J_W(x)y=W(y,x)x.$$
We say that a pseudo-Riemannian manifold is {\it conformally Osserman} if for each base point the spectrum (or equivalently the
characteristic polynomial) of the conformal Jacobi operator $J_W(x)$ is independent of the choice of the unit tangent vector $x$. Note
that the spectrum is allowed to vary from point to point. As in the case of point-wise Osserman manifolds, conformally Osserman and
self-dual manifolds of dimension $4$ are closely related. More precisely, we have the following theorem (see also
\cite{BG1},\cite{BGV}). 

\begin{theorem}\label{SD}
A $4$-dimensional oriented pseudo-Riemannian manifold is conformally Osserman if and only if it is self-dual or anti-self-dual. 
\end{theorem}

It should be pointed out that in the positive definite setting in dimension $4$ one can actually show point-wise existence of a Clifford
algebra structure compatible with a (conformal) Osserman algebraic curvature tensor \cite{BG1}. More precisely, for a conformal Osserman
algebraic  curvature tensor $R$ there exist skew-adjoint operators $\Phi_1, \Phi_2, \Phi_3=\Phi_2\Phi_1$ and constants $\ll_1,
\ll_2, \ll_3$ such that 
\begin{itemize}
\item $\Phi_i\Phi_j+\Phi_j\Phi_i=-2\delta_{ij}\id$;
\item $\ll_1+\ll_2+\ll_3=0$;
\item $W=\ll_1R_{\Phi_1}+\ll_2R_{\Phi_2}+\ll_3R_{\Phi_3}$, where
\begin{equation}\label{RPhi}
R_\Phi(x,y)z:=g(\Phi y,z)\Phi x-g(\Phi x,z)\Phi y -2g(\Phi
x,y)\Phi z.
\end{equation}
\end{itemize}
The tensors $R_\Phi$ naturally appear in the geometry of {\ekaterina rank $1$-symmetric} spaces: if we let 
$$R_0(x,y)z=g(x,z)y-g(y,z)x$$
denote the Riemann curvature tensor of the standard sphere, then the curvature tensor of the projective spaces $\mathbb{C}P^n$
and $\mathbb{H}P^n$ can be written as
$$R_0+R_I\ \ \ \ \text{and}\ \ \ \ \ R_0+R_I+R_J+R_K, \ \ \text{respectively}.$$
{\ekaterina Here $I$ and $\{I,J,K\}$ are the canonical complex and quaternion structures on $\mathbb{C}P^n$ and on
$\mathbb{H}P^n$, respectively}. Operators of this type also are central to the analysis of \cite{X4}.

In this paper we will primarily study Osserman $0$-models of signature $(2,2)$; they are a
convenient algebraic abstraction of what is happening (point-wise) on Osserman and conformal Osserman manifolds (of signature
(2,2)). Their exact definitions and other preliminaries can be found in Section \ref{prelim}. Section \ref{cliffalg} is
dedicated to the proof of the following proposition, which is our main algebraic result; {\ekaterina this result plays a crucial
role in the analysis of \cite{CGV,GGVV}.}
\begin{proposition}\label{alg}
A $0$-model $\text{\pbglie M}=(V,g,A)$ of signature $(2,2)$ is Osserman if and only if there exist skew-adjoint linear
operators $\Phi_1, \Phi_2, \Phi_3=\Phi_2\Phi_1$ and constants $\ll_i, \ll_{ij}$ such that 
\begin{itemize}
\item $\Phi_1^2=-\id, \Phi_2^2=\Phi_3^2=\id$;
\smallbreak
\item $\Phi_i\Phi_j+\Phi_j\Phi_i=0$ if $i\neq j$;
\smallbreak
\item $A=\ll_0R_0+\sum \ll_iR_{\Phi_i}+\sum_{i<j}
      \ll_{ij}[R_{\Phi_i}+R_{\Phi_j}-R_{(\Phi_i-\Phi_j)}]$.
\end{itemize}
\end{proposition}
A triple of skew-adjoint operators $(\Phi_1, \Phi_2, \Phi_3)$, $\Phi_3=\Phi_2\Phi_1$, satisfying the first two identities of
the previous proposition is often referred to as a {\it $\mathrm{Cliff}(1,1)$-structure on $(V,g)$}. A very interesting aspect
of the tensor decomposition in the third identity is the appearance of the tensor $R_{(\Phi_i-\Phi_j)}$. It is particularly
interesting that for some choices of $i,j$ we have $(\Phi_i-\Phi_j)^2=0$. Due to this distinct feature we will say that 
$A$ admits a {\it modified $\mathrm{Cliff}(1,1)$-structure}. It should also be pointed out that the usage of the tensor $R_0$
is optional (see Remark \ref{rem1} below).

Proposition \ref{alg} can be applied to the geometric setting. In Section \ref{geo} we discuss the following geometric results.
As in Proposition \ref{alg} the usage of the tensor $R_0$ is optional. 

\begin{theorem}\label{thm1}
A pseudo-Riemannian manifold of signature $(2,2)$ is point-wise Osserman if and only if for each point of the manifold there
exists a local smooth $\mathrm{Cliff}(1,1)$-structure $(\Phi_1, \Phi_2,\Phi_3)$ and smooth functions $\ll_i,\ll_{ij}$ such that 
$$R=\ll_0R_0+\sum \ll_iR_{\Phi_i}+\sum_{i<j}
      \ll_{ij}[R_{\Phi_i}+R_{\Phi_j}-R_{(\Phi_i-\Phi_j)}].$$
\end{theorem}
In the light of the Theorem \ref{SD} we have the following  characterization of the Weyl curvature of (anti-)self-dual
manifolds. 
\begin{theorem}\label{thm2}
Let $(M,g)$ be an oriented manifold of signature $(2,2)$. The following conditions are equivalent. 
 \begin{enumerate}
 \item $M$ is conformally Osserman;
 \item $M$ is self-dual  or anti-self-dual;
 \item For each $P\in M$ there exists a local smooth $\mathrm{Cliff}(1,1)$-structure $(\Phi_1, \Phi_2,\Phi_3)$ and smooth
functions $\ll_i,\ll_{ij}$ such that $\ll_0-\ll_1+\ll_2+\ll_3=0$ and 
$$W=\ll_0R_0+\sum \ll_iR_{\Phi_i}+\sum_{i<j}
      \ll_{ij}[R_{\Phi_i}+R_{\Phi_j}-R_{(\Phi_i-\Phi_j)}].$$  
\end{enumerate}
\end{theorem}

The non-homogeneous examples of (conformally) Osserman manifolds given in \cite{BGV} indicate that in general one can not find
a $\mathrm{Cliff}(1,1)$-structure $(\Phi_1, \Phi_2, \Phi_3)$ such that the corresponding functions $\ll_i, \ll_{ij}$ are
constant. Indeed, if one could find such $\Phi_i$ the Jordan normal form of the (conformal) Jacobi operators would have to be
independent of the base point, contrary to the examples of \cite{BGV}.  A natural question at this point is if the (Weyl)
curvature tensor of a {\it globally (conformally) Jordan Osserman manifold}, i.e. (conformally) Osserman manifold whose Jordan
normal form of the (conformal) Jacobi operator is independent of the base point, allows a decomposition with constant functions
$\ll_i, \ll_{ij}$. {\ekaterina An} affirmative answer to this question is proven in Section \ref{geo}.
\begin{theorem}\label{thm3}
If a connected pseudo-Riemannian manifold of signature $(2,2)$ is globally Jordan Osserman (resp. globally conformally Jordan
Osserman) then the  $\mathrm{Cliff}(1,1)$-structure of Theorem \ref{thm1} (resp. Theorem \ref{thm2})  can be chosen so that the
functions $\ll_i,\ll_{ij}$ are constant. 
\end{theorem}

\section{Preliminaries}\label{prelim}

In what follows we will assume $(M,g)$ is a $4$-dimensional pseudo-Riemannian manifold of {\ekaterina neutral} signature $(2,2)$.
When referring to a pseudo-orthonormal frame or a pseudo-orthonormal basis $\{e_1, e_2, e_3, e_4\}$ we will always assume 
$$g(e_i,e_j)=\epsilon_i\delta_{ij},\ \ \  \text{where}\ \  \epsilon_1=\epsilon_2=-1,\ \ \epsilon_3=\epsilon_4=1.$$

The Riemann curvature tensor $R$ of $(M,g)$ satisfies the following symmetries:
\begin{eqnarray}
&& R(x,y,z,v)=-R(y,x,z,v)=-R(x,y,v,z),\label{eq:csymm.a}\\
&& R(x,y,z,v)=R(z,v,x,y), \text{ and }\label{eq:csymm.b}\\
&& R(x,y,z,v)+R(y,z,x,v)+R(z,x,y,v)=0\label{eq:csymm.c}
\end{eqnarray}
The curvature tensor $R$ restricted to a tangent space $T_PM$ is an  example of an {\it algebraic curvature tensor:} a
$4$-tensor on an
{\ivacomments innerproduct space} 
which satisfies symmetries (\ref{eq:csymm.a})-(\ref{eq:csymm.c}). This
abstract setting is convenient when working with point-wise properties of the geometric curvature tensor. The triple
$\text{\pbglie M}=(V,g,A)$, where $A$ is an algebraic curvature tensor on the innerproduct space $(V,g)$, is called a $0$-{\it
model} \cite{Gi07}. 

Another important example of a $0$-model comes from conformal geometry. The {\it Weyl tensor} $W$ is obtained from the
decomposition 
\begin{equation}
R=\frac{Scal}{24}\ g\cdot g+\frac{1}{2}\Big(Ric - \frac{Scal}{4}\  g\Big)\cdot g+W,\label{defw}
\end{equation}
where $Ric$ denotes the Ricci curvature (the contraction of $R$  with respect to the first and the third slots), $Scal$ denotes the scalar curvature (the contraction of $Ric$) and where $h\cdot k$ deonotes the Kulkarni-Nomizu product of two symmetric $2$-tensors: \begin{eqnarray*}
h\cdot k\  (v_1, v_2, v_3, v_4)&=h(v_1, v_3)k(v_2, v_4)+h(v_2, v_4)k(v_1, v_3)\\
&-h(v_1, v_4)k(v_2, v_3)-h(v_2, v_3)k(v_1, v_4).
\end{eqnarray*}
The Weyl tensor depends only on the conformal class of $(M,g)$. Morover, it satisfies the curvature symmetries
(\ref{eq:csymm.a})-(\ref{eq:csymm.c}) and so we can treat it abstractly as an algebraic curvature tensor. In fact, we may use
the decomposition (\ref{defw}) to associate the Weyl tensor $W_A$ to any  $0$-model $(V,g,A)$. Note that the Weyl tensor is
always Ricci flat. 

An algebraic curvature tensor $A$ gives rise to the {\it Jacobi operator}, a family of operators $J_A(x)$ defined by 
$$g(J_A(x)y,z)=A(y,x,x,z).$$ It follows from the curvature symmetries 
(\ref{eq:csymm.a})-(\ref{eq:csymm.c}) that each $J_A(x)$, $x\neq 0$ induces a self-adjoint operator on the orthogonal
complement $\{x\}^\perp$. In particular, for unit timelike vectors $x$ in a vector space of signature $(2,2)$ the operator
$J_A(x)$ may be viewed as a self-adjoint operator on a vector space of signature $(1,2)$.  

Following the terminology of the spectral geometry of the Riemann curvature tensor, we say that a $0$-model $(V,g,A)$ is {\it
Osserman} (resp. {\it conformal Osserman}) if the characteristic polynomial of $J_A(x)$ (resp. $J_{W_A}(x)$) does not depend on
the choice of timelike unit vector $x$. If $A$ is (conformal) Osserman then the characteristic polynomial of the Jacobi operator
does not depend on the choice of unit spacelike vector $x$ either, see \cite{Gi02}. Osserman $0$-models of signature $(2,2)$
have been classified based upon the form of the corresponding minimal polynomial \cite{ABBR, BBR,BGV}. We have the following
four types. 

\begin{theorem}\label{Types}
A $0$-model $(V,g,A)$ is Osserman if and only if one of the following holds.
\begin{itemize}
\item {\sc Type I:} The Jacobi operators $J_A(x)$, $\|x\|^2=-1$, are diagonalizable, i.e. have matrix representations of the
form 
\begin{equation}
 \left[
\begin{array}{ccc}
 \alpha & 0 & 0 \\
0 & \beta & 0 \\
0 & 0 & \gamma \end{array} \right] ,\quad \alpha,\beta,\gamma\in\mathbb{R}. \label{eq:11a}
\end{equation}
In this case if $\{e_1,e_2,e_3,e_4\}$ is a pseudo-orthonormal basis with respect to which $J_A(e_1)$ has the matrix
representation as above, then the non-vanishing components of $A$ with respect to $\{e_1,e_2, e_3, e_4\}$ are:
\begin{eqnarray*}
&& A_{1221}=A_{4334}=-\alpha,\ A_{1331}=A_{4224}=\beta,\ A_{1441}=A_{3223}=\gamma,\\
&& A_{1234}=\frac{2\alpha-\beta-\gamma}{3},\ A_{1423}=\frac{-\alpha-\beta+2\gamma}{3},\ 
A_{1342}=\frac{-\alpha+2\beta-\gamma}{3}.
\end{eqnarray*}
\item {\sc Type II:} The Jacobi operators $J_A(x)$, $\|x\|^2=-1$, have matrix representations of the form
\begin{equation}
\left[
\begin{array}{ccc}
\alpha & \beta & 0 \\
-\beta & \alpha & 0  \\
 0 & 0 & \gamma
\end{array}
\right] ,\quad \alpha,\beta,\gamma\in\mathbb{R},\ \beta\neq 0.\label{eq:11b}
\end{equation}
In this case if $\{e_1,e_2,e_3,e_4\}$ is a pseudo-orthonormal basis with respect to which $J_A(e_1)$ has the matrix
representation as above, then the non-vanishing components of $A$ with respect to $\{e_1,e_2, e_3, e_4\}$ are:
\begin{eqnarray*}
&& A_{1221}=A_{4334}=-\alpha,\ A_{1331}=A_{4224}=\alpha,\ A_{1441}=A_{3223}=\gamma,\\
&& A_{2113}=A_{2443}=-\beta,\ A_{1224}=A_{1334}=\beta,\\
&& A_{1234}=\frac{\alpha-\gamma}{3},\ A_{1423}=\frac{2(\gamma-\alpha)}{3},\ A_{1342}=\frac{\alpha-\gamma}{3}.
\end{eqnarray*}
\item {\sc Type III:} The Jacobi operators $J_A(x)$, $\|x\|^2=-1$, have matrix representations of the form 
\begin{equation}
 \left[
  \begin{array}{ccc}    \epsilon ( \alpha - {\textstyle{\frac{1}{2}}} ) &
    \epsilon {\textstyle{\frac{1}{2}}} & 0 \\
    -\epsilon {\textstyle{\frac{1}{2}}} &
    \epsilon (\alpha + {\textstyle{\frac{1}{2}}}) & 0 \\
    0 & 0 & \beta
  \end{array}
\right] ,  \quad \epsilon = \pm 1,\ \alpha,\beta,\in\mathbb{R}. \label{eq:11c}
\end{equation}
In this case if $\{e_1,e_2,e_3,e_4\}$ is a pseudo-orthonormal basis with respect to which $J_A(e_1)$ has the matrix
representation as above, then the non-vanishing components of $A$ with respect to $\{e_1,e_2, e_3, e_4\}$ are:
\begin{eqnarray*}
&& A_{1221}=A_{4334}=-\epsilon(\alpha-\frac{1}{2}),\ A_{1331}=A_{4224}=\epsilon(\alpha+\frac{1}{2}),\\
&&A_{1441}=A_{3223}=\beta,\\
&& A_{2113}=A_{2443}=-\epsilon\frac{1}{2},\ A_{1224}=A_{1334}=\epsilon\frac{1}{2},\\
&& A_{1234}=\frac{\epsilon(\alpha-\frac{3}{2})-\beta}{3},\ A_{1423}=\frac{-2\epsilon\alpha+2\beta}{3},
\ A_{1342}=\frac{\epsilon(\alpha+\frac{3}{2})-\beta}{3}.
\end{eqnarray*}
\item {\sc Type IV:} The Jacobi operators $J_A(x)$, $\|x\|^2=-1$, have matrix representations of the form 
\begin{equation}
\left[
\begin{array}{ccc}
 \alpha  & 0  &
{\textstyle{\frac{\sqrt{2}}{2}}} \\
 0 & \alpha &
{\textstyle{\frac{\sqrt{2}}{2}}} \\
-{\textstyle{\frac{\sqrt{2}}{2}}} &
{\textstyle{\frac{\sqrt{2}}{2}}} & \alpha
\end{array}
\right] , \quad \alpha\in\mathbb{R}.\label{eq:11d}
\end{equation}
In this case if $\{e_1,e_2,e_3,e_4\}$ is a pseudo-orthonormal basis with respect to which $J_A(e_1)$ has the matrix
representation as above, then the non-vanishing components of $A$ with respect to $\{e_1,e_2, e_3, e_4\}$ are:
\begin{eqnarray*}
&& A_{1221}=A_{4334}=-\alpha,\ A_{1331}=A_{4224}=\alpha,\ A_{1441}=A_{3223}=\alpha,\\
&& A_{2114}=A_{2334}=-\frac{\sqrt{2}}{2},\ A_{3114}=-A_{3224}=\frac{\sqrt{2}}{2},\\
&& A_{1223}=A_{1443}=A_{1332}=-A_{1442}=\frac{\sqrt{2}}{2}.\end{eqnarray*}
\end{itemize}
\end{theorem}

Here is an important corollary of this classification result.

\begin{corollary}\label{uniqueness}
Let $A$ and $\widetilde{A}$ be two Osserman algebraic curvature tensors on an innerproduct space $(V,g)$ of signature $(2,2)$.
If for some timelike unit vector $x$ we have $J_{A}(x)=J_{\widetilde{A}}(x)$, then necessarily {\ekaterina$A=\widetilde{A}$.}
\end{corollary}

As mentioned in the Introduction, (Jordan) Osserman algebraic curvature tensors tend to be related to representations of
Clifford algebras. We proceed by investigating this relationship in signature $(2,2)$. 

\section{Clifford structures and Proposition \ref{alg}}\label{cliffalg}

Let $\big(\mathbb{R}^{(p,q)}, (.,.)\big)$ denote the standard innerproduct space of signature $(p,q)$. The Clifford algebra
$\mathrm{Cliff}(p,q)$ is the unital algebra generated by $\mathbb{R}^{(p,q)}$ subject to the Clifford commutation relations:
$$v\cdot w+w\cdot v=-2(v,w)\cdot 1.$$
{\ivacomments The Clifford algebra $\mathrm{Cliff}(0,2)$, for example,  can be seen as the algebra of quaternions  $\mathbb{H}:=\operatorname{Span}_{\mathbb{R}}\{1,i,j,k\}$, where 
$$i^2=j^2=k^2=-1,\quad ij+ji=ik+ki=jk+kj=0,\quad ijk=-1\,.$$} 
{\ekaterina Note that the
multiplication on the left by $i, j$ and $k$ gives rise to a unitary representation of $\mathrm{Cliff}(0,2)$ on $\mathbb{R}^4$. 
Likewise, $\mathrm{Cliff}(1,1)=\operatorname{Span}_{\mathbb{R}}\{I,J,K\}$ are the 
para-quaternions; these satisfy the relations
$$I^2=-1,\quad J^2=K^2=1,\quad IJ+JI=IK+KI=JK+KI=0,\quad IJK=1\,.$$}
The main results of our paper rely on the existence of the following representation of $\mathrm{Cliff}(1,1)$ on
$\mathbb{R}^{(2,2)}$.

\begin{lemma} \label{Cliff-lemma} Let $\{e_1, e_2, e_3, e_4\}$ be a pseudo-orthonormal basis for an innerproduct space $(V,g)$ of signature $(2,2)$.
There exist skew-adjoint linear maps $\Phi_i:V\to V$, where $i=1, 2, 3$, such that 
\begin{enumerate}
\item $\Phi_i\Phi_j+\Phi_j\Phi_i=0$, $i\neq j$; 
\item $\Phi_1^2=-\id$, $\Phi_2^2=\Phi_3^2=\id$; 
\item $\Phi_3=\Phi_2\Phi_1$;
\item $\Phi_1(e_1)=e_2$, $\Phi_2(e_1)=e_3$, $\Phi_3(e_1)=e_4$.
\end{enumerate}
\end{lemma}
\begin{proof} Our choice of orthonormal basis allows us to identify $V$ with $\mathbb{R}^{(1,1)}\otimes
\mathbb{R}^{(0,2)}$. More precisely, there exists an isometry $T:V\to \mathbb{R}^{(1,1)}\otimes \mathbb{R}^{(0,2)}$ with 
\begin{alignat*}{5}
&T(e_1)=(1,0)\otimes(1,0),\ \ \ T(e_2)=(1,0)\otimes(0,1),\\
&T(e_3)=(0,1)\otimes(1,0),\ \ \ T(e_4)=(0,1)\otimes(0,1).
\end{alignat*}
Consider 
\begin{equation*}
\alpha_0=\left[ \begin{array}{cc}
1 & \ 0\\
0 & -1  
\end{array} \right],\ \ \  
\alpha_1=\left[ \begin{array}{cc}
0 & 1\\
1 & 0  
 \end{array} \right],\ \ \ 
\alpha_2=\left[ \begin{array}{cc}
0 & -1\\
1 & \ 0  
 \end{array} \right]=\alpha_1\alpha_0;
\end{equation*}
these matrices satisfy $\alpha_i\alpha_j+\alpha_j\alpha_i=0$ if $i\neq j$, $\alpha_0^2=\alpha_1^2=\id$ and $\alpha_2^2=-\id$.
Viewed as operators on $\mathbb{R}^{(1,1)}$ $\alpha_0$ and $\alpha_2$ are self-adjoint while $\alpha_1$ is skew-adjoint. It is
now easy to check that $\Phi_i$, $i=1, 2, 3$, defined by 
$$T\circ\Phi_1\circ T^{-1}=\alpha_0\otimes \alpha_2,\ \ 
T\circ\Phi_2\circ T^{-1}=\alpha_1\otimes \id,\ \ 
T\circ\Phi_3\circ T^{-1}=\alpha_2\otimes \alpha_2
$$
satisfy conditions (1)-(3). {\ekaterina We note that:
$$\id(1,0)=\alpha_0(1,0)=(1,0),\ \ \ \alpha_1(1,0)=\alpha_2(1,0)=(0,1)\,.$$
Property (4) now follows.}
\end{proof}

{\ekaterina Gilkey and Ivanova \cite{GI01} gave a construction using Clifford algebras that showed the Jordan normal form of a
Jordan Osserman algebraic curvature tensor can 
{\ivacomments be }
arbitrary. We use their construction in what follows. Note that it
follows from Lemma
\ref{Cliff-lemma} that} the map
$\Phi_1$ is an isometry and that the maps
$\Phi_2,
\Phi_3$ are anti-isometries:
$$g(\Phi_2 v,\Phi_2 w)=g(\Phi_3 v,\Phi_3 w)=-g(v,w).$$
Therefore,  for unit timelike vector $x$ the set $\{x,\Phi_1 x, \Phi_2x, \Phi_3x\}$ forms a pseudo-orthonormal basis.

\begin{lemma}
Let $\Phi_i$, $i=1, 2, 3$, be skew-adjoint maps on an innerproduct space $(V,g)$ of signature $(2,2)$ satisfying relations (1)
and (2) of Lemma \ref{Cliff-lemma}. {\ekaterina Let $\ll_i$ and 
{\ivacomments $\ll_{ij}$ with  $i<j$ } 
be real constants. Then the curvature tensor 
$$R:=\lambda_0R_0+\lambda_1R_{\Phi_1}+\lambda_2R_{\Phi_2}+\lambda_3R_{\Phi_3}+
\sum_{i< j}\lambda_{ij}\big[R_{\Phi_i}+R_{\Phi_j}-R_{(\Phi_i-\Phi_j)}\big]$$
gives rise to an Osserman $0$-model on $V$.} 
\end{lemma}
\begin{proof}
In the computation which follows we will use $\pi_x$ to denote the linear map 
$$\pi_x(v):=g(v,x)x.$$
The Jacobi operator corresponding to an algebraic curvature tensor of the form $R_\Phi$ (see (\ref{RPhi})) takes the form 
$$J_{R_\Phi}(x)y=-3g(\Phi y, x)\Phi x=3g(y,\Phi x)\Phi x, \ \ \ \mathrm{i.e.}\ \ \  J_{R_\Phi}(x)=3\pi_{\Phi x}.$$
The matrix representations of the operators 
$$J_i(x):=J_{R_{\Phi_i}}(x)=3\pi_{\Phi_ix},\ \ \ J_{ij}(x):=J_{R_{(\Phi_i-\Phi_j)}}(x)=
3\pi_{(\Phi_ix-\Phi_jx)}$$
with respect to $\{x,\Phi_1x,\Phi_2x,\Phi_3x\}$ are independent of the choice of unit timelike vector $x$. 
Therefore, 
\begin{equation}\label{jacR}
J_R(x)=\ll_0(\pi_x+ \id)+ 3\sum_{i=1}^3\lambda_i\pi_{\Phi_i x}+
3\sum_{i< j}\lambda_{ij}[\pi_{{\ekaterina\Phi_ix}}+\pi_{{\ekaterina\Phi_jx}}-\pi_{(\Phi_ix-\Phi_jx)}]
\end{equation}
has its matrix representation with respect to $\{x,\Phi_1x,\Phi_2x,\Phi_3x\}$ independent of the choice of a timelike unit
vector $x$. So, the algebraic curvature tensor $R$ is Osserman.
\end{proof}

We can explicitly write down the matrix representation of the operator $J_R(x)$ (see equation (\ref{jacR})) with respect to the
basis $\{x,\Phi_1x,\Phi_2x,\Phi_3x\}$. To do so note that the operator 
$\pi_{{\ekaterina\Phi_ix}}+\pi_{{\ekaterina\Phi_jx}}-\pi_{(\Phi_ix-\Phi_jx)}$ is zero on $\mathrm{Span}\{{\ekaterina\Phi_ix},{\ekaterina\Phi_jx}\}^\perp$, while it
acts as
$${\ekaterina\Phi_ix}\mapsto \epsilon_i{\ekaterina\Phi_jx},\ \ \ \ {\ekaterina\Phi_jx}\mapsto \epsilon_j{\ekaterina\Phi_ix}$$
on  $\mathrm{Span}\{{\ekaterina\Phi_ix},{\ekaterina\Phi_jx}\}$. Hence our matrix representation is of the form 
\begin{equation}\label{matrixform}
\left[\begin{array}{ccc}
\lambda_0-3\lambda_1 & 3\lambda_{12} &3\lambda_{13}\\
-3\lambda_{12} &\lambda_0+3\lambda_2 &3\lambda_{23}\\
-3\lambda_{13} &3\lambda_{23} &\lambda_0+3\lambda_3
\end{array}\right].
\end{equation}
\bigbreak
The proof of Proposition \ref{alg} is based upon the previous two lemmas and the classification of Theorem \ref{Types}.
\medbreak
\begin{proof}[Proof of Proposition \ref{alg}]
Consider a pseudo-orthonormal basis $\{e_1, e_2, e_3, e_4\}$  for $(V,g)$ and consider maps $\Phi_i$, $i=1, 2, 3$ of  Lemma
\ref{Cliff-lemma}; we have
$$e_2=\Phi_1(e_1),\ \ e_3=\Phi_2(e_1),\ \ e_4=\Phi_3(e_1).$$
The matrix representation of the Jacobi operator $J_A(e_1)$ with respect to $\{e_2, e_3, e_4\}$ is of the form
$$J_A(e_1)=\left[\begin{array}{ccc}
\ a & b & c\\
-b & d & e\\
-c & e & f
\end{array}\right],$$
where $a=-A_{2112}$, $b=-A_{3112}$, $c=-A_{4112}$, $d=A_{3113}$, $e=A_{4113}$ and $f=A_{4114}$.  
For a suitable choice of $\ll$'s the matrix (\ref{matrixform}) reduces to matrix above. To be precise, we need
\begin{equation}\label{lambdas}
\ll_1=\frac{\ll_0-a}{3},\ \ll_2=\frac{d-\ll_0}{3},\ \ll_3=\frac{f-\ll_0}{3},\ \ll_{12}=\frac{b}{3},
\ \ll_{13}=\frac{c}{3},\ \ll_{23}=\frac{e}{3}.
\end{equation}
Let  $\lambda_i,\lambda_{ij}$ be chosen as in (\ref{lambdas}) and let 
\begin{equation}\label{CliffR}
R:=\lambda_0R_0+\lambda_1R_{\Phi_1}+\lambda_2R_{\Phi_2}+\lambda_3R_{\Phi_3}+\sum_{i< j}
\lambda_{ij}\big[R_{\Phi_i}+R_{\Phi_j}-R_{(\Phi_i-\Phi_j)}\big].
\end{equation}
It follows from the previous lemma that $R$ is Osserman with $J_R(e_1)=J_A(e_1)$. We now use Corollary \ref{uniqueness} to
conclude that $R=A$.\end{proof}

\begin{remark} \label{rem1}
The constant $\ll_0$ from the previous proof remains undetermined. Using (\ref{lambdas}) we see that 
$0=\ll_0R_0+\frac{\ll_0}{3}R_{\Phi_1}-\frac{\ll_0}{3}R_{\Phi_2}-\frac{\ll_0}{3}R_{\Phi_3}$ for all $\ll_0$ i.e. 
$$3R_0=-R_{\Phi_1}+R_{\Phi_2}+R_{\Phi_3}.$$
Moreover, we can always set $\ll_0=0$ and eliminate the $R_0$ term from the decomposition (\ref{CliffR}).
On the other hand, the modified Clifford terms $R_{\Phi_i}+R_{\Phi_j}-R_{(\Phi_i-\Phi_j)}$ are unavoidable in the cases
when the conformal Jacobi operator is non-diagonalizable.
\end{remark}
\medbreak
Proposition \ref{alg} can be applied to the Weyl tensor $W_A$ of a conformal Osserman $0$-model. Note though that the Ricci
flatness of $W_A$ imposes certain conditions on the constants $\ll_i$.  A short computation shows that the Ricci tensor
corresponding to  (\ref{CliffR}) is  $Ric=3(\ll_0-\ll_1+\ll_2+\ll_3)g$; this means we should restrict our attention to constants
$\ll_i$ with $\ll_0-\ll_1+\ll_2+\ll_3=0$.

\begin{corollary}\label{algW}
A $0$-model $\text{\pbglie M}=(V,g,A)$ of signature $(2,2)$ is conformal Osserman if and only if the Weyl tensor $W_A$ allows a
modified $\mathrm{Cliff}(1,1)$-structure 
$$W_A=\ll_0R_0+\sum \ll_iR_{\Phi_i}+\sum_{i<j}
      \ll_{ij}[R_{\Phi_i}+R_{\Phi_j}-R_{(\Phi_i-\Phi_j)}],$$
with {\ekaterina$\ll_0-\ll_1+\ll_2+\ll_3=0$.}
\end{corollary}

\section{Local geometry}\label{geo}

If the Riemann curvature tensor at a point $P\in M$ is Osserman, we say that the manifold $(M,g)$ is {\it Osserman at} $P$;
pseudo-Riemannian manifolds which are Osserman at each of their  points are called {\it point-wise Osserman}. If the Weyl
tensor of a pseudo-Riemannian manifold $(M,g)$  is Osserman then $M$ is said to be {\it conformally Osserman}. The first
examples of (point-wise) Osserman and conformal Osserman manifolds of signature $(2,2)$ are locally isotropic spaces
$\mathbb{R}^{(2,2)}$, $S^{(2,2)}$ and $\mathbb{C}P^{(1,1)}$. The para-complex projective space $\widetilde{\mathbb{C}}P^2$
\cite{Gi02} is another example of an Osserman and conformally Osserman manifold. 

The curvature decomposition of Proposition \ref{alg} carries over to the geometric setting, resulting in Theorems \ref{thm1}
and \ref{thm2}. To see this one needs to look back at the proof of Proposition \ref{alg} and replace the pseudo-orthonormal
basis $\{e_1, e_2, e_3, e_4\}$ with a smooth local pseudo-orthonormal frame. The functions $\ll_i$ and $\ll_{ij}$ in this case
are smooth because they are given in terms of the components of the curvature tensor. 

The really interesting examples of (conformally) Osserman manifolds  come from manifolds having a parallel degenerate
distribution of rank $2$: the Walker manifolds  \cite{BGV, Gi07}. These manifolds provide a wide family of (conformally)
Osserman manifolds which are not homogeneous. More specifically, we have the existence of (conformally) Osserman manifolds whose
(conformal) Jacobi operator changes its Jordan normal form from point to point. Therefore, the decomposition of the Proposition
\ref{alg} with {\it constant} $\ll$'s cannot be extended to the geometric setting without imposing further restrictions on the
(conformal) Jacobi operators. 

We say that a pseudo-Riemannian manifold is {\it globally Jordan Osserman (resp. globally conformally Jordan Osserman)} if the Jacobi operator $J_R(x)$ (resp. conformal Jacobi operator $J_W(x)$) has its Jordan normal form independent of the choice of unit timelike tangent vector $x$ and its base point. To prove Theorem \ref{thm3}, the curvature decomposition result in this setting, we need a more geometric version of the Theorem \ref{Types}.

\begin{extlemma}
Let $(M,g)$ be a globally Jordan Osserman (resp. globally conformally Jordan Osserman) manifold of signature $(2,2)$. Then for
each point of $M$ there exists a smooth local pseudo-orthonormal frame $\{e_1, e_2, e_3, e_4\}$ with respect to which the matrix
representation of the Jacobi operator $J_R(e_1)$ (resp. conformal Jacobi operator $J_W(e_1)$) is of one of the four types listed
in Theorem \ref{Types}. 
\end{extlemma}

\begin{proof} In what follows we will let $A$ denote the Riemann or the Weyl curvature tensor, depending on whether we are
working in the conformal setting. 

Consider a (smooth) unit timelike vector field $e_1$ defined on a neighborhood $\mathcal{U}_P$ of $P\in M$. For each $Q\in
\mathcal{U}_P$ there exists a pseudo-orthonormal basis $\mathcal{E}_Q=\{e_2, e_3, e_4\}$ of ${\ekaterina\{e_1\}^\perp\subset
T_QM}$ such that the matrix representation $L$ of $J_A(e_1)$ with respect to $\mathcal{E}_Q$ is independent of $Q$ and is of one
of the four types listed in the Theorem \ref{Types}. A priori we do not know if  $\mathcal{E}_Q$ is smooth in $Q$. 

Let $\mathcal{F}:=\{x_2, x_3, x_4\}$, with $x_2$ timelike,  be a smooth pseudo-orthonormal frame for the sub-bundle
${\ekaterina\{e_1\}^\perp\subset TM\big|_{\mathcal{U}_P}}$. Let $X_Q$ be the change of basis matrix at $Q\in\mathcal{U}_P$:
$$[e_2, e_3, e_4]=[x_2, x_3, x_4]\cdot X_Q.$$
As in the case of $\mathcal{E}_Q$ it is not clear if $X_Q$ is smooth in $Q$.  Let $L_Q$ be the matrix representation of the
Jacobi operator $J_A(e_1)$ with respect to the frame $\mathcal{F}$ at a point $Q\in\mathcal{U}_P$; the entries of $L_Q$ depend
smoothly on $Q$.  Note that 
\begin{equation}\label{system}
X_QL=L_QX_Q \ \ \ \ \text{and}\ \ \ \ \ X_Q^TGX_Q=G,
\end{equation}
where $G=\mathrm{diag}(-1, 1, 1)$. Without loss of generality we will  assume $\mathcal{E}_P=\mathcal{F}$, i.e. $X_P=\id$ and
$L=L_P$. 

To prove our Extension Lemma it suffices to show the existence of a solution $X_Q$ of the system (\ref{system}) which is smooth
in $Q$ (on a neighborhood of $P$) and which satisfies $X_P=\id$.

We start with the first of our two equations. Consider the family of linear transformations
$$T_Q:Y\mapsto YL-L_QY,\ \ Q\in\mathcal{U}_P$$
on the vector space of $3\times 3$ matrices. The dimension of $\ker{T_Q}$ is independent of $Q$ due to 
$$T_Q(Y)=YL-X_QLX_Q^{-1}Y=X_Q\big(X_Q^{-1}YL-LX_Q^{-1}Y\big)=X_Q\cdot T_P (X_Q^{-1} Y)$$
and the fact that multiplications by invertable matrices are linear isomorphisms. Set $$k:=\dim \ker{T_Q}=\dim \ker
(\mathrm{ad} L),$$ where $\mathrm{ad}L=[L,.]$. 
By {\ekaterina Cramer's} Rule the $9\times 9$ system of equations 
\begin{equation}
X_Q\cdot L-L_Q\cdot X_Q=0
\end{equation}
has a $k$-parameter family of solutions $X_Q=X_Q(\vec{\ll})$. Since the coefficients of the system vary smoothly with $Q$ the
solutions $X_Q(\vec{\ll})$ depend smoothly on $Q$ and linearly on $\vec{\ll}$. Note that we can always find parameters
$\vec{\ll_0}$ so that 
{\ivacomments $X_P(\vec{\ll_0})=\id$}.

Our next step is to use the second equation of (\ref{system}) to solve for (some of) the parameters $\vec{\ll}$ in the form of
smooth functions of $Q$. We will accomplish this via the Implicit Function Theorem.
 
Let  $\mathrm{Symm}_G=\{S\ |\ GS^TG=S\}$. Consider  the function $$F:\mathcal{U}_P\times \mathbb{R}^k\to \ker (\mathrm{ad}L)
\cap \mathrm{Symm}_G$$ given by
$$F:(Q,\vec{\ll})\mapsto GX_Q(\vec{\ll})^TGX_Q(\vec{\ll}).$$
This function is well-defined, i.e. $\mathrm{Im}\ F\subseteq \ker (\mathrm{ad}L) \cap \mathrm{Symm}_G$, due to  
a straightforward computation involving identities $LG=GL^T$ and $L_QG=GL_Q^T$.

Recall that $X_Q(\vec{\lambda})$ depends linearly on $\vec{\lambda}$ and that $X_P(\vec{\lambda}_0)=\id$. This means that the
linearization of $F$ with respect to $\vec{\ll}$ at $(P,\vec{\ll_0})$ is 
\begin{equation}\label{lin}
\mathcal{L}:\vec{\ll}\mapsto GX_P(\vec{\ll})^TG+X_P(\vec{\ll}).
\end{equation}
By the Implicit Function Theorem it suffices to prove that the map (\ref{lin}) is onto. To understand this map's rank and
nullity note that $\vec{\ll}\mapsto X_P(\vec{\ll})$ is an isomorphism between $\mathbb{R}^k$ and 
{\ivacomments $\ker(\mathrm{ad}L)$.}
With this in
mind it is clear that 
$$\dim \ker\mathcal{L}=\dim \Big(\ker(\mathrm{ad}L)\cap \text{{\pbglie{so}}}_G\Big),$$
where $\text{{\pbglie{so}}}_G=\{ S\ |\ GS^TG=-S \}.$

A short computation involving $LG=GL^T$ shows that if $X\in \ker (\mathrm{ad}L)$ then also  $GX^TG\in \ker (\mathrm{ad}L)$.
Since
$X=\frac12(X-GX^TG)+\frac12(X+GX^TG)$
we have that
$$\ker(\mathrm{ad}L)=\Big(\ker(\mathrm{ad}L)\cap \text{{\pbglie{so}}}_G\Big)\oplus\Big(\ker (\mathrm{ad}L) \cap
\mathrm{Symm}_G\Big).$$ Consequently,
\begin{equation*}
k=\dim \ker(\mathrm{ad}L)=\dim \Big(\ker(\mathrm{ad}L)\cap \text{{\pbglie{so}}}_G\Big)+\dim\Big(\ker (\mathrm{ad}L) \cap
\mathrm{Symm}_G\Big)
\end{equation*}
and 
{\ivacomments $k-\dim\ker\mathcal{L}=\dim \big(\ker (\mathrm{ad}L) \cap \mathrm{Symm}_G\big)$. }
Since $\mathrm{Im}\ F\subseteq  \ker (\mathrm{ad}L)
\cap \mathrm{Symm}_G$ the linear map (\ref{lin}) is onto. 
{\ivacomments Applying }
the Implicit Function Theorem we obtain (some of) the
parameters $\vec{\ll}$ as smooth functions of $Q$; this yields to a 
{\ivacomments local }
smooth solution $X_Q$ of (\ref{system}). The proof of our
Extension Lemma is now complete.
\end{proof}

The proof of the last of our results, Theorem \ref{thm3}, is immediate from the Extension Lemma and the proof of Proposition
\ref{alg}. One simply needs to replace the pseudo-orthonormal basis $\{e_1, e_2, e_3, e_4\}$ with the smooth local
pseudo-orthonormal frame given by the Extension Lemma.

\section{Related problems}
We conclude our paper with two related open problems.
\begin{enumerate}
\item Modified Clifford algebraic curvature tensor is always (both spacelike and timelike) Jordan Osserman. It is known that in
certain dimensions and signatures (such as signature (2,2) presented here) the converse also holds: a (Jordan) Osserman
algebraic curvature tensor allows a modified Clifford structure. Is this a phenomenon which holds in general?
\medbreak
\item Corollary \ref{uniqueness} is the crucial background result which made our approach successful. Is there a higher
dimensional analogue of Corollary \ref{uniqueness}? Such a result would prove useful in answering the question raised above. 
\end{enumerate}

\section*{Acknowledgments} Research of P. Gilkey partially supported by Project MTM2006-01432 (Spain). Research of S. Nik\v
cevi\'c partially supported by Project 144032 (Srbija).


\begin{thebibliography}{AAA}



\bibitem{ABBR} D. Alekseevsky, N. Bla{\v z}i{\'c}, N. Bokan, Z. Raki{\'c},
{\em  Self-duality and pointwise Osserman spaces}, Arch. Math. (Brno),
{\bf 35}, (1999), 193--201.

\bibitem{BG1}N. Bla{\v z}i{\'c}, P. Gilkey, {\it Conformally
Osserman manifolds and self-duality in Riemannian geometry}, Differential geometry and its applications, Proceedings , 15--18, MATFYZPRESS, Prague, 2005.



 \bibitem{BBR} N. Bla\v zi\'c, N. Bokan and Z. Raki\'c,
  {\it Osserman pseudo-Riemannian manifolds of signature $(2,2)$}, J.Austral. Math. Soc.
 {\bf 71} (2001), 367--395.

\bibitem{X3} N. Bla\v zi\'c and P. Gilkey, {\it Conformally Osserman manifolds and conformally complex space forms},
Int. J. Geom. Method. in Modern Physics {\bf 1} (2004), 97--106.

\bibitem{X7} J.P. Bourguignon, {\it Some highlights of Dmitri Alekseevsky's work}, 
Int. J. Geom. Method. in Modern Physics {\bf 3} (2006), 823--831.


\bibitem{BGV} M. Brozos-V\'{a}zquez, E. Garc\'{\i}a--R\'{\i}o, R.
V\'{a}zquez-Lorenzo,  {\it Conformally Osserman four-dimensional
manifolds whose conformal Jacobi operators have complex eigenvalues}, 
Proc. R. Soc. Lond. Ser. A Math. Phys. Eng. Sci. {\bf 462} (2006), no. 2069, 1425--1441

\bibitem{CGV} E. Calvi\~{n}o-Louzao, E. Garc\'{\i}a-R\'{\i}o, and R. V\'{a}zquez-Lorenzo,
{\it Four-dimensional Osserman--Ivanov--Petrova metrics of neutral signature},
{Class. Quantum Grav.} {\bf 24} (2007), 2343--2355.


\bibitem{C} Q. S. Chi, 
{\it A curvature characterization of certain locally rank one symmetric spaces},
{J. Differential Geom.} {\bf 28}
(1988), 187--202.

\bibitem{X8} L. Cordero and P. Parker, {\it Lattices and periodic geodesics in pseudo-Riemannian 2-step nilpotent
Lie groups}, Int. J. Geom. Method. in Modern Physics {\bf 5} (2008), 79--99.

\bibitem{X4} J. D\'ias-Ramos, B. Fiedler, E. Garc\'{\i}a-R\'{\i}o, and P. Gilkey,
{\it The structure of algebraic covariant derivative curvature tensors}, 
Int. J. Geom. Method. in Modern Physics {\bf 1} (2004), 711--720.

\bibitem{X1} E. Garc\'{\i}a-R\'{\i}o, A. Haji-Badali, M. E. Vazquez-Abal, R. Vazquez-Lorenzo,
{\it Lorentzian 3-manifolds with commuting curvature operators}, Int. J. Geom. Method. in Modern Physics {\bf 5} (2008), 557-572.

\bibitem{GGVV} E. Garc\'{\i}a-R\'{\i}o, P. Gilkey, M. E. V\'{a}zquez-Abal, and R.
V\'{a}zquez-Lorenzo,
{\it Four-dimensional Osserman metrics of neutral signature, arXiv:0804.0436}.

\bibitem{Gi02} P. Gilkey, {\bf Geometric properties of natural operators
defined by the Riemann curvature tensor}, World Scientific
Publishing Co., Inc., River Edge, NJ, 2001.

\bibitem{Gi07} P. Gilkey, {\bf The Geometry of Curvature Homogeneous Pseudo-Riemannian Manifolds}, 
Advanced Texts in Mathematics, 2. Imperial College Press, London, 2007.

\bibitem{GI01}{\ekaterina P. Gilkey and R. Ivanova,  {\it The Jordan normal form of Osserman algebraic curvature tensors},
Results in Math. {\bf 40} (2001) 192--204.}

\bibitem{X5} P. Gilkey and S. Nik\v cevi\'c, {\it Affine curvature homogeneous $3$-dimensional Lorentz manifolds},
Int. J. Geom. Method. in Modern Physics {\bf 2} (2005), 737--749.

\bibitem{X2} P. Gilkey and S. Nik\v cevi\'c, {\it Pseudo-Riemannian Jacobi-Videv Manifolds},
Int. J. Geom. Method. in Modern Physics {\bf 4} (2007), 727--738.

 \bibitem{GSV} P. Gilkey, A. Swan and L. Vanhecke, {\it Isoparametric  geodesic spheres and a
 conjecture of Osserman regarding the  Jacobi Operator}
 Quart. J. Math. Oxford 46 (1995), 299--320.


\bibitem{Nik0} Y. Nikolayevsky, {\it Osserman manifolds and Clifford structures}, Houston J. Math. {\bf 29}
(2003), no.1, 59--75 (electronic).

\bibitem{Nik1} Y. Nikolayevsky,
{\it Two theorems on Osserman manifolds}, Differential Geom. Appl.
{\bf 18} (2003), 239--253.

\bibitem{Nik2} Y. Nikolayevsky, {\it Osserman Conjecture in dimension $n \ne 8, 16$}, Math. Ann. {\bf 331} (2005), no. 3, 505--522.

\bibitem{Oss} R. Osserman,
  {\it Curvature in the eighties},
  Amer. Math. Monthly, {\bf97}, (1990) 731--756.


\bibitem{SV} K.  Sekigawa,  and  L. Vanhecke, {\it  Volume
preserving geodesic symmetries on four
   dimensional  K\"ahler manifolds},
   Differential geometry  Pen\'scola, 1985, Proceedings
   (A. M. Naveira, A. Fernandez and F. Mascaro,eds.),
    Lecture Notes in Math. {\bf  1209}, Springer,  275--290.

\end{thebibliography}
\end{document}